\documentclass[12pt]{amsart}
\usepackage[T1]{fontenc}
\usepackage[utf8]{inputenc}
\usepackage{lmodern}
\usepackage{amsmath,amssymb}
\usepackage{tikz}
\usetikzlibrary{positioning,calc}
\usepackage[margin=1in]{geometry}
\usepackage{microtype}
\usepackage{enumitem}
\usepackage[hidelinks]{hyperref}

\newtheorem{theorem}{Theorem}[section]
\newtheorem{proposition}[theorem]{Proposition}
\newtheorem{lemma}[theorem]{Lemma}

\newtheorem{conjecture}[theorem]{Conjecture}
\theoremstyle{definition}
\newtheorem{definition}[theorem]{Definition}
\newtheorem{example}[theorem]{Example}
\newtheorem{remark}[theorem]{Remark}

\DeclareMathOperator{\reg}{reg}
\DeclareMathOperator{\Tor}{Tor}
\newcommand{\kk}{k}
\newcommand{\NN}{\mathbb N}
\newcommand{\ZZ}{\mathbb Z}
\newcommand{\RR}{\mathcal R}
\newcommand{\GG}{\mathcal G}
\newcommand{\K}{\mathrm K}
\newcommand{\HH}{\mathrm H}
\newcommand{\xx}{\mathbf x}
\newcommand{\cc}{\mathbf c}
\newcommand{\zz}{\mathbf z}
\newcommand{\alphavec}{\boldsymbol\alpha}

\numberwithin{equation}{section}

\newcommand{\eqmark}[1]{\refstepcounter{equation}\label{#1}}
\newcommand{\eqnumber}{\eqno{\textup{(\theequation)}}}

\setlist[enumerate]{label=\textup{(\roman*)},leftmargin=2.1em,itemsep=2pt,topsep=4pt}
\title[The Nevo--Peeva conjecture]{The Nevo--Peeva conjecture on powers of edge ideals}
\author{T\`ai Huy H\`a}
\address{Mathematics Department, Tulane University, New Orleans, LA 70118, USA}
\email{tha@tulane.edu}
\date{}
\hypersetup{pdftitle={The Nevo--Peeva conjecture on powers of edge ideals},pdfauthor={Tài Huy Hà}}

\begin{document}
\begin{abstract}
We establish the conjecture of Nevo and Peeva, which asserts that all sufficiently large powers of the edge ideal of a gapfree graph have linear resolutions. Equivalently, gapfreeness characterizes the graphs for which some positive power of the edge ideal has a linear resolution. 
\end{abstract}
\maketitle
\raggedbottom


\section{Introduction}\label{sec:intro}

Let $G$ be a finite \emph{simple} graph (i.e., $G$ has neither loops nor multiple edges) with vertex set $V(G)=\{x_1,\ldots,x_n\}$ and edge set $E(G)$, and let $S=\kk[x_1,\ldots,x_n]$ be a polynomial ring over a field $\kk$. The edge ideal of $G$ is 
$$I(G)=(x_ix_j\mid \{x_i,x_j\}\in E(G)).$$
Throughout the paper, we assume that $G$ has at least one edge. We are interested in the minimal graded free resolutions of the powers of $I(G)$, and in how the combinatorial structure of $G$ controls the degrees of their syzygies.

Since $I(G)^q$ is generated in degree $2q$, it has a linear resolution if and only if $\reg I(G)^q=2q$. The asymptotic results of Cutkosky, Herzog and Trung~\cite{CHT} and Kodiyalam~\cite{Kod} imply that $\reg I(G)^q=2q+b$ for all sufficiently large $q$, where $b$ is a nonnegative integer. Finding graph theoretic conditions for $b=0$ is a basic problem in the study of powers of edge ideals; see~\cite{BKHsurvey} for background and related questions.

A graph is called \emph{gapfree} if it has no induced subgraph consisting of two disjoint edges. An observation of Francisco, H\`a and Van Tuyl, recorded in~\cite[Proposition~1.8]{NP}, shows that gapfreeness is necessary for any positive power of $I(G)$ to have a linear resolution. Nevo and Peeva~\cite{NP} constructed a gapfree graph whose edge ideal has a square without a linear resolution. Thus, gapfreeness does not force linear resolutions to begin with the second power. They raised the following conjecture for sufficiently large powers.

\begin{conjecture}[Nevo--Peeva {\cite[Open Problem~1.11(1)]{NP}}]\label{conj:NP}
Let $G$ be a finite simple gapfree graph with at least one edge, and let $I=I(G)$ be its edge ideal over a field. Then $I^q$ has a $2q$-linear resolution for all sufficiently large integers $q$.
\end{conjecture}

The purpose of this paper is to prove Conjecture~\ref{conj:NP}. Together with the necessary condition above, our result gives the following characterization.

\begin{theorem}\label{thm:main}
Let $G$ be a finite simple graph with at least one edge, and let $I=I(G)$ be its edge ideal over a field $\kk$. The following conditions are equivalent.
\begin{enumerate}
\item $G$ is gapfree.
\item There exists an integer $q_0\geq1$ such that $I^q$ has a $2q$-linear resolution for every $q\geq q_0$.
\item $I^q$ has a linear resolution for some integer $q\geq1$.
\end{enumerate}
\end{theorem}

The implication \textup{(i)}$\Rightarrow$\textup{(ii)} is the main assertion. In terms of asymptotic regularity, it says that the constant $b$ vanishes exactly when $G$ is gapfree.

Several important cases of Conjecture~\ref{conj:NP} have been established. Fr\"oberg's theorem~\cite{Froberg} characterizes edge ideals with linear resolutions by chordality of the complementary graph. Herzog, Hibi and Zheng~\cite{HHZ} showed that, in this case, all powers have linear resolutions. For gapfree graphs, Bigdeli, Herzog and Zaare-Nahandi~\cite{BHZ} proved that all powers are linearly presented. The remaining focus is therefore the behavior of higher nonlinear syzygies.

Banerjee~\cite{Banerjee} introduced the concept of \emph{even-connections} to study colon ideals arising from powers and proved linearity for powers of gapfree graphs without induced \emph{crickets}. Erey~\cite{Erey18,Erey19} established corresponding results under the exclusion of an induced \emph{diamond} or a cycle of length four. Banerjee, Kara and H\`a~\cite{BKH} developed an approach using local regularity, proving that $\reg I(G)^q=2q$ for all $q\geq2$ when $G$ is gapfree and $\reg(I(G):x)\leq2$ for every vertex $x$. Minh and Vu~\cite{MV} subsequently characterized linear resolutions of the second and third powers in terms of gapfreeness and the regularity of the original edge ideal.

Another approach is to investigate the stronger property of \emph{linear quotients}. Erey, Faridi, H\`a, Hibi, Kara and Morey~\cite{EFHHKM} studied operations preserving this property, constructed families whose powers have linear quotients, and proposed a persistence conjecture for linear quotients of powers. More recently, Erey, H\`a, Hibi, Kara and Morey~\cite{EHHKM} proved that every power $I(G)^q$, $q\geq2$, has linear quotients when $G$ is gapfree with matching number three. 

Our method is to establish the linear quotients property for auxiliary ideals associated to affine charts of the Rees algebra of $I(G)$.
The main combinatorial ingredient is Theorem~\ref{thm:matching}. We consider a graph $H$ and a family $L$ of label sets on its vertices. Each unmatched vertex is assigned a label from its label set, and each matched pair contributes $uv$. If every label set contains $u$ or $v$, the resulting ideal $J(H,L)$ has linear quotients in decreasing lexicographic order. The proof combines an augmenting path argument with an interpolation argument for the exponents of $u$ and $v$.

To apply this result, fix an edge $uv$ of a gapfree graph $G$ and set $C=S[I(G)/(uv)]$. Let $W$ consist of the vertices adjacent to neither endpoint, and put $X=V(G)\setminus(\{u,v\}\cup W)$. The set $W$ is independent. For each exponent vector $\cc\in\NN^X$, we construct a graph $H_{\cc}$ with label family $L_{\cc}$ and configuration ideal $J_{\cc}=J(H_{\cc},L_{\cc})$ in $R=\kk[u,v,w\mid w\in W]$, and prove
$$
C_{\cc}=\frac{\xx^{\cc}}{(uv)^{|\cc|}}J_{\cc},
$$
where $\xx^{\cc}=\prod_{x\in X}x^{c_x}$ and $|\cc|=\sum_{x\in X}c_x$. Here $C_{\cc}$ is the component determined by the exponents on $X$. The linear resolution of $J_{\cc}$, together with the grading shift in this identity, gives a vanishing theorem for the Koszul homology of $C$ on the vertex variables.

Applying this calculation on every edge chart shows that the nonlinear part of the Koszul homology of the Rees algebra vanishes after each edge variable is inverted. Since this part is finitely generated over the polynomial ring in the edge variables, it has finite length. Its components in sufficiently large power degree therefore vanish, proving Theorem~\ref{thm:main}. The use of homogeneous localizations is related to the approach to asymptotic regularity in the author's earlier work~\cite{Ha}; here the necessary vanishing follows directly from the linear quotients theorem and the chart description.

The paper is outlined as follows. In Section~\ref{sec:matching}, we give the definition and prove the linear quotients theorem for configuration ideals. The main result of this section is Theorem~\ref{thm:matching}. In Section~\ref{sec:charts}, we identify these ideals with the components of affine charts, compute the relevant Koszul homology, and prove our main result, Theorem~\ref{thm:main}, establishing Conjecture~\ref{conj:NP}.

\medskip
\noindent\textbf{Acknowledgments.} The author learned a great deal about linear quotients during an AIM SQuaRE program with Nursel Erey, Sara Faridi, Takayuki Hibi, Selvi Kara, and Susan Morey. The idea of using linear quotients for the configuration ideals in the present paper grew out of our earlier joint work \cite{EFHHKM,EHHKM}. The author thanks his collaborators for many stimulating discussions, and the American Institute of Mathematics for its support and hospitality. The author also thanks Trung Chau, Arvind Kumar, A.V. Jayanthan and Aniketh Sivakumar for informing him of their progress and the very recent preprint of Dupont García on the same conjecture. The author was partially supported by a Simons Foundation grant.

\medskip
\noindent\textbf{Declaration of AI-use.} A generative AI was used as an assisted writing tool during the preparation of this paper, particularly, for language editing and reference cross-check. 


\section{Matchings and linear quotients}\label{sec:matching}

We prove a linear quotients theorem for monomial ideals defined by matchings and assignments of labels to unmatched vertices. These ideals will describe the components of the affine charts in Section~\ref{sec:charts}. We begin with the notation and criterion for linear quotients used in the proof.

Recall that a nonzero finitely generated graded module $M$ over a standard graded polynomial ring $R$ has a \emph{$d$-linear resolution} if 
$$\Tor_i^R(M,\kk)_j=0 \text{ for } j\ne d+i.$$
Its Castelnuovo--Mumford \emph{regularity} is 
$$\reg M=\max\{j-i\mid\Tor_i^R(M,\kk)_j\ne0\}.$$
In particular, a module generated in degree $d$ has a linear resolution exactly when its regularity is $d$.

Let $J$ be a monomial ideal in a polynomial ring $R$. We denote its minimal monomial generating set by $\GG(J)$. For monomials $M,N$, write $M:N=M/\gcd(M,N)$, as in~\cite{EFHHKM}. For a variable $z$, let $\deg_z(M)$ denote its exponent in $M$. An ordering $g_1>\cdots>g_t$ of $\GG(J)$ is a \emph{linear quotients order} if $(g_1,\ldots,g_{\ell-1}):g_\ell$ is generated by variables for every $2\leq \ell\leq t$. We recall the following well-known criterion; see, for example, \cite[Lemma~2.1]{EFHHKM}.

\begin{lemma}\label{lem:lq-criterion}
Fix an ordering $g_1>\cdots>g_t$ of the minimal monomial generators of an ideal $J$. This ordering gives linear quotients if and only if, for every $a<b$, there is an index $c<b$ such that $g_c:g_b$ is a variable dividing $g_a:g_b$. If, in addition, all generators have degree $d$, then $J$ has a $d$-linear resolution.
\end{lemma}

We now introduce the configuration ideals that will be used. Recall that a \emph{matching} of a graph is a set of pairwise disjoint edges. A matching \emph{saturates} a set of vertices if every vertex in that set belongs to one of its edges.

Let $s$ be a nonnegative integer, let $w_1,\ldots,w_s,u,v$ be distinct indeterminates, and set $R=\kk[w_1,\ldots,w_s,u,v]$. Let $H$ be a simple graph on a set $\Omega$ of $d$ labeled vertices, and let $L=(L(a))_{a\in\Omega}$ be a family of subsets of $\{w_1,\ldots,w_s,u,v\}$. We call $L(a)$ the \emph{label set} of $a$ and assume that
\eqmark{eq:endpoint-condition}
$$
L(a)\cap\{u,v\}\ne\varnothing\quad\text{for every }a\in\Omega.
\eqnumber
$$

We associate a monomial ideal to the graph $H$ together with its label sets $L$.

\begin{definition}\label{def:configuration}
A \emph{configuration} for $(H,L)$ is a pair $(M,\ell)$, where $M$ is a matching of $H$ and $\ell$ assigns to each vertex $a\in\Omega\setminus V(M)$ a label $\ell(a)\in L(a)$. Its monomial is
$$
m(M,\ell)=(uv)^{|M|}\prod_{a\in\Omega\setminus V(M)}\ell(a).
$$
Let $\GG$ be the set of distinct configuration monomials. The \emph{configuration ideal} of $(H,L)$ is $J(H,L)=(\GG)\subseteq R$.
\end{definition}

A configuration with $a$ singletons and $b$ pairs has $d=a+2b$ and produces a monomial of degree $a+2b=d$. The condition~\eqref{eq:endpoint-condition} guarantees that configurations exist, and distinct configuration monomials have the same degree. Thus $\GG=\GG(J(H,L))$. For $d=0$, we take $\GG=\{1\}$ and $J(H,L)=R$.

\begin{example}\label{ex:pairs}
	Let $\Omega=\{a,b\}$, let $L(a)=L(b)=\{u\}$, and let $ab\in E(H)$. The configuration ideal is $J(H,L)=(u^2,uv)$. The generator $uv$ comes from the pair $\{a,b\}$ and would be lost if one considered only singleton assignments. 
\end{example}

\begin{theorem}\label{thm:matching}
With the notation above, $J(H,L)$ has linear quotients in decreasing lexicographic order with respect to $w_1>\cdots>w_s>u>v$. Consequently, $J(H,L)$ has a $d$-linear resolution over every field.
\end{theorem}

We first prove two lemmas. The first treats exchanges at the labels $w_i$. The second shows that, after their exponents are fixed, the possible exponents at $u$ form an interval.

\begin{lemma}\label{lem:W-exchange}
Fix an integer $j$ with $1 \le j \le s$. Let $A,B\in\GG$. Suppose that $\deg_{w_i}(A)=\deg_{w_i}(B)$ for every $1\le i<j$, and that $\deg_{w_j}(A)>\deg_{w_j}(B)$. Then there is a variable $z$ later than $w_j$ in the order $w_1>\cdots>w_s>u>v$ such that $z\mid B$ and $w_jB/z\in\GG$.
\end{lemma}

\begin{proof}
Let $b_i=\deg_{w_i}(B)$ for $1\leq i\leq j$. For each $i$ introduce $b_i$ distinct slots labeled $w_i$, taking no slots when $b_i=0$, and let $Q$ be their union. Join a slot labeled $w_i$ to $a\in\Omega$ precisely when $w_i\in L(a)$. The configuration representing $B$ determines a matching $M_B$ saturating $Q$. Introduce one additional slot $q_*$ labeled $w_j$, and put $Q'=Q\cup\{q_*\}$. Since the exponents of $A$ and $B$ agree at $w_1,\ldots,w_{j-1}$ and $\deg_{w_j}(A)\geq b_j+1$, a configuration representing $A$ gives a matching $M_A$ saturating $Q'$.

Consider the component $P$ containing $q_*$ in the symmetric difference $M_A\triangle M_B$. The new slot is unmatched by $M_B$ and matched by $M_A$, whereas every slot in $Q$ is matched by both matchings. Thus $P$ is an alternating path ending at a vertex $a\in\Omega$. Since the endpoints of $P$ lie in opposite parts of the bipartite graph, $P$ has odd length. Its first and last edges belong to $M_A$, so $a$ is unmatched by $M_B$. Set
$$
M'=(M_B\setminus E(P))\cup(M_A\cap E(P)).
$$
Every vertex of $Q$ remains matched, $q_*$ becomes matched, and every element of $\Omega$ previously covered by $M_B$ remains covered. The only new element covered is $a$. This is the path flip shown in Figure~\ref{fig:augmenting}.

\begin{figure}[htbp]
\centering
\begin{tikzpicture}[x=1.16cm,y=1cm,
    point/.style={circle,draw,fill=white,inner sep=2pt,minimum size=7mm},
    original/.style={line width=0.85pt},
    alternative/.style={line width=0.95pt,dashed}]
\node[point] (q0) at (0,0.75) {$q_*$};
\node[point] (a1) at (1,-0.65) {$a_1$};
\node[point] (q1) at (2,0.75) {$q_1$};
\node[point] (a2) at (3,-0.65) {$a_2$};
\node[point] (q2) at (4,0.75) {$q_2$};
\node[point] (a3) at (5,-0.65) {$a_3$};
\draw[alternative] (q0)--(a1);
\draw[original] (a1)--(q1);
\draw[alternative] (q1)--(a2);
\draw[original] (a2)--(q2);
\draw[alternative] (q2)--(a3);
\node[font=\small,anchor=east] at (-0.38,0.75) {$Q'$};
\node[font=\small,anchor=east] at (-0.38,-0.65) {$\Omega$};
\draw[original] (0.2,-1.55)--(0.8,-1.55);
\node[font=\small,anchor=west] at (0.9,-1.55) {edge of $M_B$};
\draw[alternative] (3.15,-1.55)--(3.75,-1.55);
\node[font=\small,anchor=west] at (3.85,-1.55) {edge of $M_A$};
\end{tikzpicture}
\caption{The alternating path from the new point $q_*$ to an element $a_3$ not used by $M_B$ in Lemma~\ref{lem:W-exchange}. Flipping replaces the solid edges by the dashed edges. The other matching edges are unchanged.}
\label{fig:augmenting}
\end{figure}
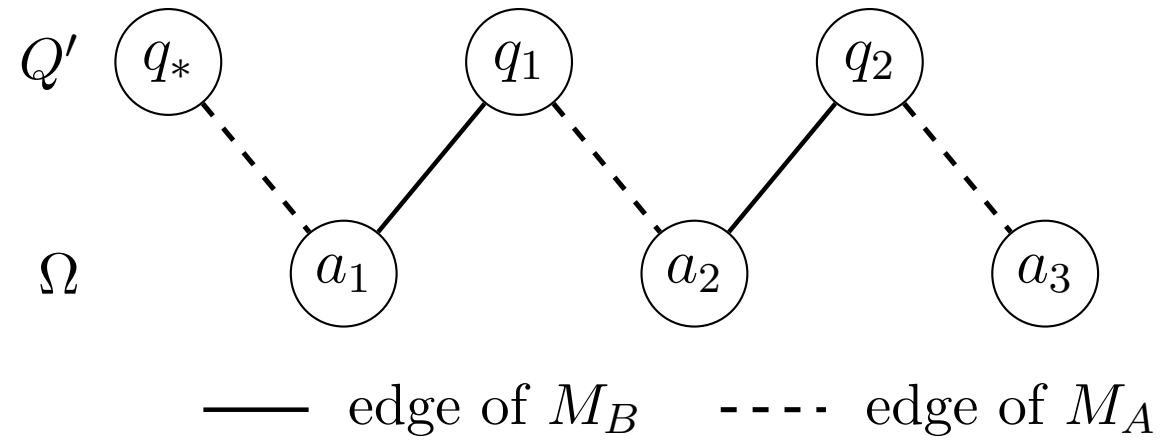

Use $M'$ in place of $M_B$ for all assignments to slots of $Q'$, leaving the assignments outside $Q$ unchanged except as specified below. The slot labels in $Q'$ differ from those in $Q$ by exactly one occurrence of $w_j$. Thus this replacement preserves the exponents at $w_1,\ldots,w_{j-1}$ and increases the exponent at $w_j$ by one, although the individual elements assigned to the old slots may change. In particular, the terminal element $a$ need not be matched to the new slot $q_*$.

Since $a$ was not assigned to a slot in $Q$, either it was a singleton with label $z\in\{w_{j+1},\ldots,w_s,u,v\}$, or it belonged to a pair. In the singleton case, remove its former label $z$; replacing the assignments on $Q$ by $M'$ then gives the monomial $w_jB/z$.

In the pair case, write the pair as $\{a,b\}\in E(H)$. Neither member was covered by $M_B$, and only $a$ is added to the covered elements by $M'$, so $b$ remains free. Delete the pair and give $b$ a label $t\in L(b)\cap\{u,v\}$, which exists by~\eqref{eq:endpoint-condition}. Setting $z=uv/t\in\{u,v\}$, the resulting monomial is
$$
\frac{B}{uv}\,w_jt=\frac{w_jB}{z}.
$$
In both cases, all copies of $\Omega$ are used exactly once, all assignments are permitted, and $z$ comes later than $w_j$. Thus $w_jB/z\in\GG$, as required.
\end{proof}

\begin{lemma}\label{lem:interval}
Fix $\alphavec=(\alpha_1,\ldots,\alpha_s)\in\NN^s$ for which a configuration monomial has exponent $\alpha_i$ at each $w_i$. Then the set
$$
\{\deg_u(M)\mid M\in\GG,\ \deg_{w_i}(M)=\alpha_i\text{ for }1\leq i\leq s\}
$$
is an interval of integers.
\end{lemma}

\begin{proof}
Fix configurations $\mathcal C$ and $\mathcal D$ having the prescribed exponents $\alpha_i$ at $w_i$. We connect them by configurations with the same $w_i$ exponents, changing the exponent of $u$ by at most one at each step.

First consider a matched pair $\{a,b\}$ in either configuration. Replace it by singletons labeled $t_a\in L(a)\cap\{u,v\}$ and $t_b\in L(b)\cap\{u,v\}$. This replaces $uv$ by $t_at_b$ in the configuration monomial, and
$$
\deg_u(t_at_b)-\deg_u(uv)\in\{-1,0,1\}.
$$
The exponents at the $w_i$ do not change. Repeating this operation yields configurations $\mathcal C'$ and $\mathcal D'$ consisting only of singletons. It is enough to connect $\mathcal C'$ to $\mathcal D'$ and then reverse the chosen unpairing steps for $\mathcal D$.

For $1\leq i\leq s$, take $\alpha_i$ distinct slots labeled $w_i$, and let $Q$ be their union. The assignments to $Q$ in $\mathcal C'$ and $\mathcal D'$ define matchings $M$ and $N$, both saturating $Q$, in the bipartite graph of slots and elements of $\Omega$. Every component of $M\triangle N$ is an alternating cycle or path. The components are disjoint, so they may be flipped successively without changing the others. A cycle flip changes only the assignments to the slots; in particular, it leaves the configuration monomial unchanged.

For an alternating path $P$, orient its endpoints as $a,b\in\Omega$, with $a$ covered by the current matching but not by $N$, and $b$ covered by $N$ but not by the current matching. Replace the edges of the current matching on $P$ by the edges of $N$ on $P$. Each slot of $Q$ remains matched, while the set of covered elements of $\Omega$ changes by replacing $a$ with $b$. Before the flip, $b$ is a singleton with some endpoint label $t_b\in L(b)\cap\{u,v\}$. Remove its singleton label and assign the released element $a$ a label $t_a\in L(a)\cap\{u,v\}$. In terms of the configuration monomial $m$, this exchange gives
$$
m\longmapsto m\frac{t_a}{t_b},
\qquad \deg_u(t_a)-\deg_u(t_b)\in\{-1,0,1\}.
$$
Thus each path flip preserves every $w_i$ exponent and changes the exponent of $u$ by at most one. This uses the same alternating path operation as Figure~\ref{fig:augmenting}, but here $Q$ is fixed and the flip exchanges one element of $\Omega$ for another.

After every component is flipped, the matching to $Q$ is exactly $N$. All remaining elements are singletons with endpoint labels. Change their labels, one element at a time, to those in $\mathcal D'$. Each change is permitted because both the current and target labels lie in the corresponding label set, and it changes the $u$ exponent by at most one. We have now reached $\mathcal D'$, so the pair removals originally chosen for $\mathcal D$ can be reversed, restoring precisely its original pairs and labels. We obtain the claimed sequence from $\mathcal C$ to $\mathcal D$.

The exponents of $u$ along this sequence are integers whose consecutive differences have absolute value at most one. They therefore take every integer between the exponents of the initial and final configurations. Since $\mathcal C$ and $\mathcal D$ were arbitrary, the set in the statement is an interval.
\end{proof}

\begin{proof}[{\bf Proof of Theorem~\ref{thm:matching}}]
Write $J=J(H,L)$. The case $d=0$ is immediate, so assume $d>0$. Write the generators in decreasing lexicographic order as $g_1>\cdots>g_t$. Fix $A=g_a$ and $B=g_b$ with $a<b$, and let $z$ be the first variable at which their exponents differ. Then $\deg_z(A)>\deg_z(B)$.

If $z=w_j$, Lemma~\ref{lem:W-exchange} gives a later variable $z'$ dividing $B$ such that $D=zB/z'\in\GG$. If $z=u$, the exponents at all $w_i$ agree. Lemma~\ref{lem:interval} gives a generator with these same exponents and with $u$ exponent $\deg_u(B)+1$. Every generator has degree $d$, so this generator is $D=uB/v$; in particular, $v\mid B$. The first difference cannot occur at $v$, since equality of the preceding exponents and equality of total degrees would force equality there as well.

In either case, $D$ precedes $B$ and $D:B=z$ divides $A:B$. Lemma~\ref{lem:lq-criterion} proves that the given ordering has linear quotients and that $J$ has a $d$-linear resolution.
\end{proof}

\begin{remark}\label{rem:pairs}
The graph $H$ is arbitrary in Theorem~\ref{thm:matching}. In particular, the interpolation argument does not require an exchange property for its matchings. 
\end{remark}


\section{Rees algebras and eventual linear resolutions}\label{sec:charts}

We apply Theorem~\ref{thm:matching} to the affine charts of the Rees algebra of an edge ideal. The resulting Koszul homology calculation will imply Theorem~\ref{thm:main} by finite generation.

\subsection{Affine charts and Koszul homology}\label{subsec:charts}

We identify a vertex with its corresponding variable and write $xy\in E(G)$ for $\{x,y\}\in E(G)$. The \emph{open neighborhood} of $x$ is $N_G(x)=\{y\in V(G)\mid xy\in E(G)\}$. A set of vertices is \emph{independent} if it contains no edge.

For a graded module $M$, our \emph{shift} convention is $M(a)_j=M_{j+a}$. If $\zz$ is a sequence of degree one elements in a graded algebra $C$, we write $\K(\zz;C)$ for its \emph{Koszul complex} and $\HH_i(\zz;C)$ for its homology. Each Koszul basis vector has homological degree one and degree one in the grading of $C$. If $R$ is the polynomial ring on the entries of $\zz$ and $M$ is a graded $R$-module, then 
$$\HH_i(\zz;M)=\Tor_i^R(M,\kk).$$
We refer the reader to~\cite{HH} for further background on monomial ideals and their resolutions.

Fix an edge $uv$ of $G$, and set
$$
W=V(G)\setminus\bigl(\{u,v\}\cup N_G(u)\cup N_G(v)\bigr),
\quad X=V(G)\setminus(\{u,v\}\cup W).
$$
In this subsection, we assume that $W$ is independent. This holds for every edge of a gapfree graph, since an edge inside $W$ would form a gap with $uv$. Set $R=\kk[u,v,w\mid w\in W]$ and consider the chart
$$
C=S[I(G)/(uv)]=S\bigl[ab/(uv)\mid ab\in E(G)\bigr]
\subseteq S[u^{-1},v^{-1}].
$$
Equip $C$ with the $\ZZ$-grading by total degree inherited from the Laurent polynomial ring $S[u^{-1},v^{-1}]$, in which every vertex variable has degree one and every edge ratio has degree zero. The degree of a Laurent monomial is therefore the sum of its exponents. This grading is nonnegative because the displayed generators of $C$ have degrees zero and one.

Every vertex in $X$ is adjacent to $u$ or $v$. Moreover, every edge other than $uv$ has at least one endpoint in $X$: there are no edges inside $W$ or between $W$ and $\{u,v\}$. All exponents on $X$ are nonnegative, so they give $C$ an additional multigrading by $\NN^X$. Every element of $R$ has $X$ multidegree zero; thus multiplication by $R$ preserves each multihomogeneous component $C_{\cc}$. Since the total grading and the $X$ multigrading are compatible, each $C_{\cc}$ inherits a $\ZZ$-grading from the total degree, and the $R$-action respects this grading. Consequently, $C=\bigoplus_{\cc\in\NN^X}C_{\cc}$ as $\ZZ$-graded $R$-modules.

Fix $\cc=(c_x)_{x\in X}$ and write $d=|\cc|=\sum_{x\in X}c_x$ and $\xx^{\cc}=\prod_{x\in X}x^{c_x}$. Let $\Omega_{\cc}$ contain $c_x$ distinct copies of each $x\in X$, and let $\pi:\Omega_{\cc}\to X$ send each copy to its original vertex. Define the graph $H_{\cc}$ and its label sets $L_{\cc}=(L_{\cc}(a))_{a\in\Omega_{\cc}}$ by
$$
E(H_{\cc})=\bigl\{\{a,b\}\subseteq\Omega_{\cc}\mid \pi(a)\pi(b)\in E(G)\bigr\},
\quad L_{\cc}(a)=N_G(\pi(a))\cap(\{u,v\}\cup W).
$$
In particular, distinct copies of the same vertex are never adjacent in $H_{\cc}$. Set $J_{\cc}=J(H_{\cc},L_{\cc})\subseteq R$. Since $\pi(a)\in X$, the set $L_{\cc}(a)$ contains $u$ or $v$. Thus Theorem~\ref{thm:matching} implies that $J_{\cc}$ has a $d$-linear resolution.

\begin{figure}[htbp]
\centering
\begin{tikzpicture}[x=1cm,y=0.95cm,
    vert/.style={circle,draw,fill=white,minimum size=7mm,inner sep=1pt},
    lab/.style={font=\small}]
\node[lab] at (1.15,2.75) {$G$};
\node[vert] (u) at (0.25,1.80) {$u$};
\node[vert] (v) at (2.05,1.80) {$v$};
\node[vert] (x) at (0.35,0.55) {$x$};
\node[vert] (y) at (2.05,0.55) {$y$};
\node[vert] (w) at (0.35,-0.75) {$w$};
\draw[line width=1.4pt] (u)--(v);
\draw (u)--(x) (u)--(y) (x)--(y) (x)--(w);
\node[lab,anchor=north] at (1.15,-1.16) {$W=\{w\},\quad X=\{x,y\}$};
\draw[->,line width=0.8pt] (3.15,0.64)--(4.0,0.64);
\node[lab] at (5.45,2.75) {$H_{\cc},\ \cc=(1,1)$};
\node[vert] (xc) at (4.75,1.35) {$x^{(1)}$};
\node[vert] (yc) at (6.25,1.35) {$y^{(1)}$};
\draw[line width=0.9pt] (xc)--(yc);
\node[lab,anchor=north] at (5.4,0.63) {$L_{\cc}(x^{(1)})=\{u,w\}$};
\node[lab,anchor=north] at (5.4,0.05) {$L_{\cc}(y^{(1)})=\{u\}$};
\node[lab,anchor=north] at (5.4,-0.65) {$J_{\cc}=(u^2,uw,uv)$};
\end{tikzpicture}
\caption{The graph $G$, the graph $H_{\cc}$ of labeled copies, and the configuration ideal for $\cc=(1,1)$. The pair $\{x^{(1)},y^{(1)}\}$ contributes $uv$; the singleton assignments contribute $u^2$ and $uw$. The component $C_{\cc}$ is obtained by multiplying $J_{\cc}$ by $xy/(uv)^2$.}
\label{fig:chart-construction}
\end{figure}
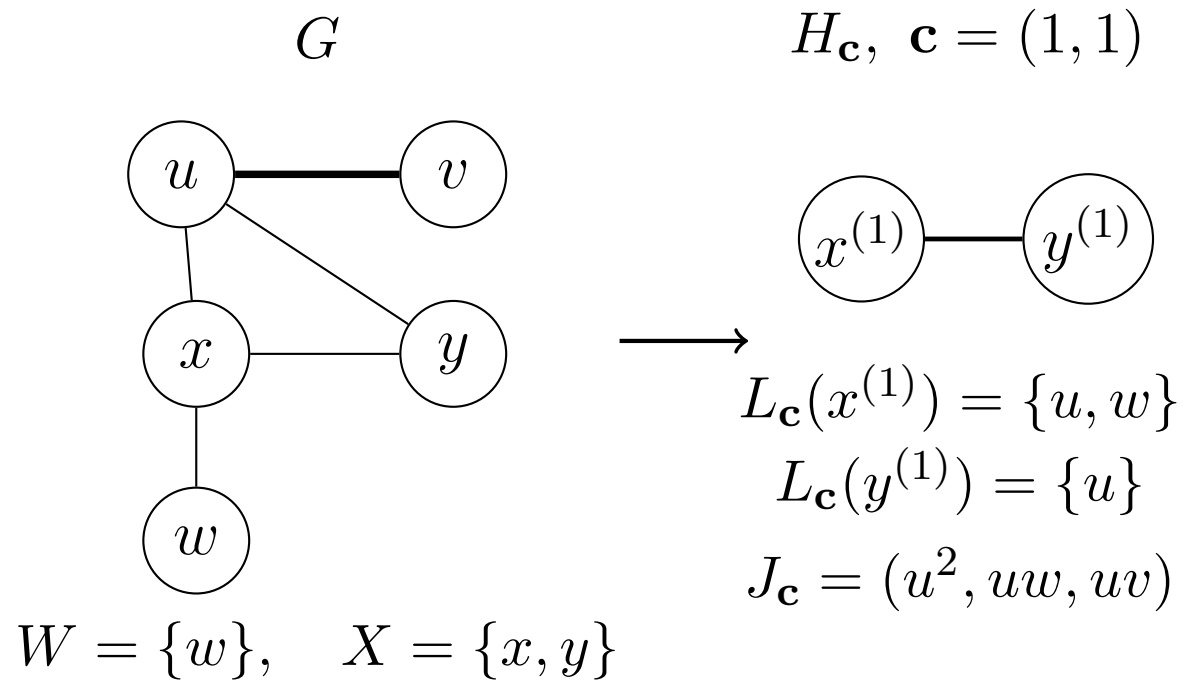

The construction is illustrated in Figure~\ref{fig:chart-construction}. In this example, the edge $xy$ inside $X$ permits a pair, whose contribution to a configuration monomial is $uv$.

\begin{proposition}\label{prop:chart}
For every $\cc\in\NN^X$, there is an equality of $R$-submodules of $S[u^{-1},v^{-1}]$,
\eqmark{eq:chart}
$$
C_{\cc}=\frac{\xx^{\cc}}{(uv)^{|\cc|}}J_{\cc}.
\eqnumber
$$
Consequently, multiplication by $\xx^{\cc}/(uv)^{|\cc|}$ gives an isomorphism $J_{\cc}(|\cc|)\cong C_{\cc}$ of $\ZZ$-graded $R$-modules, where the grading on $C_{\cc}$ is the total grading.
\end{proposition}

\begin{proof}
Write $d=|\cc|$. We prove both inclusions of~\eqref{eq:chart}. Let $(M,\ell)$ be a configuration for $(H_{\cc},L_{\cc})$, and put $b=|M|$. There are $d-2b$ singleton vertices. Each pair $\{a,a'\}\in M$ gives the edge ratio $\pi(a)\pi(a')/(uv)$, and each singleton $a$ contributes $\pi(a)\ell(a)/(uv)$. All these fractions belong to $C$ by the definitions of $H_{\cc}$ and $L_{\cc}$. Since every copy of $\Omega_{\cc}$ occurs exactly once, their product is
$$
\begin{aligned}
\prod_{\{a,a'\}\in M}\frac{\pi(a)\pi(a')}{uv}
\prod_{a\in\Omega_{\cc}\setminus V(M)}\frac{\pi(a)\ell(a)}{uv}
&=\frac{\xx^{\cc}}{(uv)^{d-b}}
\prod_{a\in\Omega_{\cc}\setminus V(M)}\ell(a)\\
&=\frac{\xx^{\cc}}{(uv)^d}(uv)^b
\prod_{a\in\Omega_{\cc}\setminus V(M)}\ell(a)
=\frac{\xx^{\cc}}{(uv)^d}m(M,\ell).
\end{aligned}
$$
The multiplier on the left has $X$ multidegree $\cc$. Multiplication by any monomial of $R$ preserves this multidegree, proving $\xx^{\cc}(uv)^{-d}J_{\cc}\subseteq C_{\cc}$.

Conversely, it suffices to consider a monomial $f\in C_{\cc}$. Express $f$ as a product of ordinary vertex variables and edge ratios, omitting any factors $uv/(uv)=1$. Every edge other than $uv$ has at least one endpoint in $X$. For an ordinary factor $x\in X$, choose $p\in\{u,v\}$ adjacent to $x$, and let $p'$ be the other endpoint. The identity
$$
x=p'\frac{px}{uv}
$$
replaces this factor by an element of $R$ times an edge ratio. Applying this substitution to every ordinary factor from $X$ gives
$$
f=r\prod_{k=1}^{h}\frac{e_k}{uv},
\qquad r\text{ a monomial of }R,
$$
where every numerator edge $e_k\ne uv$ meets $X$. Because the denominator involves only $u,v$, the total number of occurrences of each $x\in X$ among the $e_k$ is exactly $c_x$.

Assign these occurrences bijectively to the distinct copies of $x$ in $\Omega_{\cc}$, for each $x\in X$. An edge $e_k=xy$ with $x,y\in X$ gives the pair of copies assigned to its two occurrences. By definition, this pair is an edge of $H_{\cc}$. An edge $e_k=xz$ with $x\in X$ and $z\in\{u,v\}\cup W$ gives the singleton copy assigned to $x$, labeled $z\in L_{\cc}(a)$. Since every occurrence is used exactly once, these pairs and singletons form a configuration $(M,\ell)$ of $(H_{\cc},L_{\cc})$. The displayed product calculation now gives
$$
f=r\frac{\xx^{\cc}}{(uv)^d}m(M,\ell)
\in\frac{\xx^{\cc}}{(uv)^d}J_{\cc},
$$
which proves the reverse inclusion.

Finally, $\xx^{\cc}/(uv)^d$ is a nonzero Laurent monomial of total degree $d-2d=-d$. Multiplication by it is injective in $S[u^{-1},v^{-1}]$ and sends degree $j$ of $J_{\cc}$ to total degree $j-d$ of $C_{\cc}$. Under our shift convention, it therefore gives the stated isomorphism $J_{\cc}(d)\cong C_{\cc}$ of $\ZZ$-graded $R$-modules.
\end{proof}

We next compute Koszul homology on this chart. Let $\zz$ be the sequence consisting of $u,v$ and the variables in $W$, and let $\mathbf v$ be the sequence of all vertex variables of $G$.

\begin{theorem}\label{thm:chart-koszul}
Assume that the vertices adjacent to neither endpoint of $uv$ form an independent set. For $C=S[I(G)/(uv)]$, one has
\eqmark{eq:chart-vanishing}
$$
\HH_i(\zz;C)_r=\HH_i(\mathbf v;C)_r=0
\quad\text{whenever }r\ne i.
\eqnumber
$$
\end{theorem}

\begin{proof}
Multiplication by the entries of $\zz$ preserves the $X$ grading. Thus the decomposition $C=\bigoplus_{\cc\in\NN^X}C_{\cc}$ gives a direct sum decomposition of $\K(\zz;C)$. Homology commutes with direct sums, and Proposition~\ref{prop:chart} yields
$$
\HH_i(\zz;C)\cong
\bigoplus_{\cc\in\NN^X}\Tor_i^R(J_{\cc},\kk)(|\cc|).
$$
By Theorem~\ref{thm:matching}, the summand indexed by $\cc$ is concentrated in total degree $i$ after this shift. This proves the assertion for $\zz$.

For each $x\in X$, choose $a_x\in\{u,v\}$ adjacent to $x$, and let $b_x$ be the other endpoint. Then $x=\lambda_x b_x$, where $\lambda_x=a_xx/(uv)\in C_0$ and $C_0$ denotes the component of total degree zero. In the free module defining $\K(\mathbf v;C)$, replace the basis vector $e_x$ by $e_x-\lambda_xe_{b_x}$ for every $x\in X$, and keep the basis vectors indexed by $u,v,W$ unchanged. This is an invertible triangular change of basis over $C$. It is homogeneous for the $\ZZ$-grading, and the new vectors indexed by $X$ have differential zero.

It follows that $\K(\mathbf v;C)$ is isomorphic to the tensor product of $\K(\zz;C)$ and the exterior complex on $|X|$ zero entries. Each exterior generator has homological and total degree one. Therefore
$$
\HH_i(\mathbf v;C)_r\cong
\bigoplus_{j=0}^{\min\{i,|X|\}}
\HH_{i-j}(\zz;C)_{r-j}^{\,\oplus\binom{|X|}{j}}.
$$
The first part of the proof shows that a summand can be nonzero only when $r-j=i-j$, or $r=i$.

The direct sum calculation above uses the multigrading by the exponents on $X$. The vanishing in (\ref{eq:chart-vanishing}), however, considers the total degree, which is preserved by the change of basis.
\end{proof}

\begin{example}\label{ex:C5}
Let $G$ be the cycle with edges $uv,va,aw,wb,bu$, and choose the edge $uv$. Then $W=\{w\}$ and $X=\{a,b\}$. There is no edge between $a$ and $b$, so no pairs are allowed. For $c_a=i$ and $c_b=j$, the ideals and chart components are
$$
J_{\cc}=(v,w)^i(u,w)^j \text{ and } C_{\cc}=\frac{a^ib^j}{(uv)^{i+j}}(v,w)^i(u,w)^j.
$$
Although $J_{\cc}$ is generated in degree $i+j$, the minimum total degree of $C_{\cc}$ is zero. The example also distinguishes the chart calculation from a statement about $I(G)$ itself: the edge ideal of this cycle does not have a linear resolution, even though Theorem~\ref{thm:chart-koszul} applies on every edge chart.
\end{example}

\subsection{Proof of the main theorem}\label{sec:rees}

We now pass from the chart calculation to the powers of the edge ideal. Set $I=I(G)$ and $\RR=S[It]$, and give $\RR$ the bigrading determined by $\deg x_z=(1,0)$ and $\deg t=(-2,1)$. Thus each edge generator $et$ has bidegree $(0,1)$ and
\eqmark{eq:rees-piece}
$$
\RR_{(a,q)}=(I^q)_{a+2q}t^q\quad(q\geq0).
\eqnumber
$$
In particular, $\RR$ is nonnegatively graded in both coordinates. For an element of bidegree $(a,q)$, we call $a$ its \emph{normalized degree} and $q$ its \emph{power degree}. Thus the normalized degree is the ordinary polynomial degree minus $2q$.

Let $Q=\kk[T_e\mid e\in E(G)]$, where $\deg T_e=(0,1)$. The homomorphism 
$$S[T_e\mid e\in E(G)]\longrightarrow\RR$$ sending $T_e$ to $et$ is surjective. Write $H_i=\HH_i(\mathbf v;\RR)$ for the Koszul homology on all vertex variables, with every Koszul basis vector having bidegree $(1,0)$.

\begin{lemma}\label{lem:finite}
For every $i$, $H_i$ is a finitely generated bigraded $Q$-module. Its components satisfy
\eqmark{eq:Tor-slices}
$$
(H_i)_{(a,q)}\cong\Tor_i^S(I^q,\kk)_{a+2q}
\quad(q\geq0).
\eqnumber
$$
Moreover, $(H_i)_{(a,q)}=0$ for $a<i$.
\end{lemma}

\begin{proof}
The Koszul complex is a finite complex of finitely generated modules over the noetherian ring $S[T_e\mid e\in E(G)]$. Its homology is therefore finitely generated over this ring. Multiplication by each vertex variable is null-homotopic, with homotopy given by exterior multiplication by the corresponding Koszul basis vector. Hence every vertex variable annihilates $H_i$, and $H_i$ is finitely generated over the quotient $Q$.

By~\eqref{eq:rees-piece}, the component of the Koszul complex in power degree $q$ is the usual Koszul complex of $I^q$, with total degrees reduced by $2q$. This proves~\eqref{eq:Tor-slices}. Finally, $\RR$ has no negative normalized degrees, so the $i$th Koszul module, and therefore its homology, has no normalized degree below $i$.
\end{proof}

\begin{lemma}\label{lem:localization}
Let $e=uv\in E(G)$ and put $C_e=S[I/e]$. There is a bigraded isomorphism
\eqmark{eq:localized-Rees}
$$
\RR[(et)^{-1}]\cong C_e[Z,Z^{-1}],\quad Z=et,
\eqnumber
$$
where $\deg Z=(0,1)$ and the total grading of $C_e$ agrees with the normalized degree. If the vertices adjacent to neither $u$ nor $v$ form an independent set, then $(H_i)_{T_e}$ is concentrated in normalized degree $i$.
\end{lemma}

\begin{proof}
The localized Rees algebra contains $S$, $Z^{\pm1}$, and every ratio $f/e=(ft)/(et)$ for $f\in E(G)$. Conversely, $ft=(f/e)Z$, so these elements generate the localization. The element $Z=et$ is algebraically independent over $C_e\subseteq S[e^{-1}]$: a Laurent polynomial relation in $Z$ would give a relation in the indeterminate $t$ over the fraction field of $S$. This proves~\eqref{eq:localized-Rees}, with the stated grading.

Localization commutes with Koszul homology. Since adjoining an invertible indeterminate is flat, the displayed isomorphism identifies $(H_i)_{T_e}$ with
$$
\HH_i(\mathbf v;C_e)\otimes_{C_e}C_e[Z,Z^{-1}].
$$
Under the independence assumption, Theorem~\ref{thm:chart-koszul} places the first factor in total degree $i$. Powers of $Z$ do not change this total degree, proving the last assertion.
\end{proof}

The next proposition supplies a single threshold for all nonlinear components, rather than a threshold depending on the degree.

\begin{proposition}\label{prop:uniform}
Assume that $G$ is gapfree. Then there is an integer $q_0\geq1$ such that $\Tor_i^S(I^q,\kk)_j=0$ for every $q\geq q_0$, every $i\geq0$, and every $j>2q+i$.
\end{proposition}

\begin{proof}
For each $i$, let
\eqmark{eq:nonlinear}
$$
N_i=\bigoplus_{a>i}\ \bigoplus_{q\geq0}(H_i)_{(a,q)}.
\eqnumber
$$
The edge variables have normalized degree zero, so $N_i$ is a bigraded $Q$-submodule of $H_i$. Lemma~\ref{lem:finite} and the noetherian property of $Q$ imply that $N_i$ is finitely generated. Since $G$ is gapfree, Lemma~\ref{lem:localization} applies to every edge and gives $(N_i)_{T_e}=0$ for all $e\in E(G)$.

Suppose $N_i\ne0$ and choose homogeneous generators $h_1,\ldots,h_m$. For each edge $e$, localization shows that a power of $T_e$ kills each $h_j$. Taking the maximum of these finitely many exponents gives an integer $a_e\geq1$ with $T_e^{a_e}N_i=0$. If $D_i=1+\sum_{e\in E(G)}(a_e-1)$, then every monomial in the $T_e$ of degree at least $D_i$ is divisible by some $T_e^{a_e}$. Therefore 
$$Q_+^{D_i}N_i=0,$$
where $Q_+=(T_e\mid e\in E(G))$.

Let $b_i$ be the largest power degree of the chosen generators. A homogeneous element of $N_i$ of power degree $q\geq b_i+D_i$ is a sum of terms $f_jh_j$, with each $f_j$ homogeneous of degree at least $D_i$ in the edge variables. Each such term is zero. Thus every component of $N_i$ in sufficiently large power degree vanishes, uniformly in the normalized degree. If $N_i=0$, no bound is needed.

Only the homological degrees $0\leq i\leq n$ occur. Choose $q_0\geq1$ exceeding all the thresholds just obtained. The interpretation~\eqref{eq:Tor-slices} now gives the assertion.
\end{proof}

We are now ready to prove our main result.

\begin{proof}[{\bf Proof of Theorem~\ref{thm:main}}]
Suppose $G$ is gapfree. By Proposition~\ref{prop:uniform}, for $q\geq q_0$ one has $\Tor_i^S(I^q,\kk)_j=0$ when $j>2q+i$. Lemma~\ref{lem:finite} gives the corresponding vanishing for $j<2q+i$. Thus $I^q$ has a $2q$-linear resolution for all $q\geq q_0$, proving (i)$\Rightarrow$(ii). The implication (ii)$\Rightarrow$(iii) is immediate.

For completeness, we include the standard gap obstruction for the converse. Suppose that the induced subgraph on $\{u,v,a,b\}$ consists of the two edges $uv$ and $ab$, and fix $q\geq1$. Set $L=(uv)^qab$. The only minimal generators of $I^q$ dividing $L$ are $g_0=(uv)^q$ and $g_1=(uv)^{q-1}ab$. Indeed, every edge appearing in such a divisor must be supported on these four vertices, and the exponents of $a,b$ allow at most one copy of $ab$.

The relation $(ab)e_{g_0}-(uv)e_{g_1}$ is a first syzygy in multidegree $L$. It is minimal: the least common multiple of $g_0,g_1$ is $L$, and no proper divisor of $L$ is divisible by two minimal generators of $I^q$. Hence no syzygy in a strictly smaller multidegree can generate this relation. Its total degree is $2q+2$, rather than $2q+1$, so $I^q$ does not have a linear resolution. Since $q$ was arbitrary, (iii) implies (i).
\end{proof}

\begin{remark}\label{rem:scope}
Theorem~\ref{thm:main} gives eventual linear resolutions, not eventual linear quotients of the powers. The latter property is proved here for the auxiliary ideals $J_{\cc}$, and does not follow for $I^q$ from the chart calculation. It would be desirable to know if $I^q$ has linear quotients for $q \gg 0$. Likewise, the proof does not give an explicit graph theoretic bound on $q_0$. Such a bound would require control of the degrees of generators and annihilating powers for the modules $N_i$ in~\eqref{eq:nonlinear}. 
\end{remark}

\end{document}